\documentclass[12pt]{amsart}
                   
\usepackage{amsmath}
\usepackage{amssymb}
\usepackage{latexsym}

                     \newcommand{\Cc}{{\mathbb{C}}}

                     \newcommand{\Dd}{{\mathbb{D}}} 
                      
                     \newcommand{\Zz}{\mathbb{Z}} 
                     \newcommand{\al}{\alpha}

\newcommand{\cala}{{\mathcal{A}}} 
 
\newcommand{\cald}{{\mathcal{D}}}

\newcommand{\call}{{\mathcal{L}}} 
 
\newcommand{\calo}{{\mathcal{O}}} 
 
\newcommand{\cals}{{\mathcal{S}}}

\newcommand{\eit}{e^{i \theta}}

                     \newtheorem{theorem}{Theorem} 
                     \newtheorem{lemma}[theorem]{Lemma} 
                     \newtheorem{corollary}[theorem]{Corollary} 
                     \newtheorem{proposition}[theorem]{Proposition}

\begin{document} 

\title{\bf Algebras generated by two bounded
 holomorphic functions} 

\author{Michael I. Stessin  \ \ \ Pascal J. Thomas
} 

\begin{abstract}
We study the closure in the Hardy space or the disk algebra of algebras
generated by two bounded functions, of which one is a finite Blaschke 
product.
We give necessary and sufficient conditions for density or finite 
codimension (of
the closure) of such algebras.  The conditions are expressed in terms of 
the inner
part of some function which is explicitly derived from each pair of 
generators.  Our
results are based on identifying $z$-invariant subspaces included in the 
closure of
the algebra.

 \end{abstract}

\maketitle


\footnotetext{The second-named author would like to thank
the University at Albany, Harvard University, and 
Brown University for their hospitality
during the completion of this work.} 

\section{Introduction}
\label{intro}
Let $g$ be an inner function in the unit disk. It is possible to
prove that
if $h\in H^2$ satisfies the condition $h(z_1)=h(z_2)$ for every
pair $z_1,z_2\in \Dd $ such that $g(z_1)=g(z_2)$, then $h$
is in the $H^2$-closure of the algebra of polynomials in $g$, \ $\Cc[g]$,
and vice versa (the proof rests on the fact that, in the 
case where $g(0)=0$, the self-mapping $g$ of the 
circle is Lebesgue measure-preserving). In fact, this result holds for a 
class of functions much broader
than inner functions. It means that
the $H^2$-closure of $\Cc[g]$ is determined by the family of level sets
of the function $g$. In particular, if $g$ is a Riemann mapping of $\Dd $
onto a polynomially convex domain $\Omega $ with a reasonable boundary, 
the algebra $\Cc[g]$ is dense
in $H^2 (\Dd)$ (on this topic, the interest reader may consult 
\cite{KLS}). 

It is natural to ask if a similar result holds for algebras generated by
two bounded functions. In particular, if a pair of functions $g_1, \ g_2$
separates points of the unit disk, under what condition is the algebra 
$\Cc[g_1,g_2]$ dense in $H^2$? 

The case of closures of algebras generated by two disk-algebra functions 
(closure
is taken in the disk-algebra metric) was intensively studied (cf 
\cite{JoW}, \cite{RGB}, \cite{SiWe}). Some conditions
for the density of such an algebra were found (we will discuss some of those 
results in comparison with
those obtained in the present paper in section 7 below).

In a special case when both $g_1$ and $g_2$ are polynomials the result
similar to the one for algebras with one generator holds.   Namely, the 
following Theorem follows directly 
from Proposition 1 of \cite{JLS}.

\vspace{.4cm}

{\bf Theorem A.}  {\em Let $p_1$ and $p_2$ be two polynomials such that 
$|p_1^\prime (z)|+|p_2^\prime (z)|>0$
for all $z\in \Cc $ and their degrees are mutually prime. Then the algebra 
$\overline{\Cc[p_1,p_2]}$ 
has finite codimension in $H^2$. If, in addition, $p_1$ and $p_2$ separate 
points of the unit disk,
then $\Cc[p_1,p_2]$ is dense in $H^2$.}

\vspace{.4cm}

The conditions that the derivatives of $p_1$ and $p_2$ do not vanish 
simultaneously
and the degrees of $p_1$ and $p_2$ are mutually prime are too restrictive. 
As we will show below, 
they could be replaced with one much less restrictive condition.

If the generators are not polynomials, the condition of point separation 
is not sufficient
for the density of the algebra in $H^2$. The following simple example 
illustrates this.
Let $\chi (z)=\exp \{ \frac{z+1}{z-1} \}$ be the singular inner function 
determined by the
unit point mass at $1$, and $\psi (z)=z\chi (z)$. Clearly, $\chi $ and 
$\psi $ separate points
of the unit disk. At the same time, every polynomial in $\chi $ and $\psi 
$ belongs to
the subspace $L=\{ 1\} + \chi H^2$, where $\{ 1\}$ is the one-dimensional 
subspace of $H^2$
consisting of constants, and $\chi H^2$ is the $z$-invariant subspace 
generated by $\chi $.
Since the codimension of $\overline{\chi H^2}$ is infinite, this implies 
that
$\Cc[\chi ,\psi ] $ is not dense in $H^2$. 

In this example both generators have non-trivial singular parts. What 
happens if one of them 
is analytic in the closed disk? The following argument pertains to the 
case when one of the generators 
is a Blaschke product of order two. It shows that even in this case the 
$H^2$ 
codimension of the closure of
an algebra could be infinite, while its generators separate points of the unit 
disk.
The study of this simple case also gives the flavor of the work which we 
carry out
in the body of the paper, where we deal with the case where one generator 
is a Blaschke product 
of any finite order $n \ge 2$.

Let $|a|<1$ and $B(z)=cz \varphi (z)$, where $\varphi 
(z)=\frac{a-z}{1-\bar{a}z}$ is a M\"obius 
involution of the disk, and $c$ a constant of modulus 
$1$. Given any $z 
\in \Dd$, there is at most another point $w$ of the 
disk such that 
$B(w)=B(z)$, given by $w =\varphi(z)$. There is 
exactly one $z_0 \in \Dd$ 
such that $z_0 =\varphi(z_0)$, and it is the only 
zero of $B'$. In the  simplest case where $B(z)=z^2$, $\varphi(z)=-z$ and 
$z_0=0$. 

 Any function $h \in H^p(\Dd) $ which satisfies 
$h\circ \varphi (z)= h(z)$ can be rewritten $h=h_1 \circ B$, where $h_1$ 
is well-defined from the set-theoretic point of view because of the 
hypothesis, is holomorphic  outside of the critical point of $B$ by taking local 
inverses, and is   actually holomorphic  everywhere by 
removing the isolated  singularity. If $p=\infty$, $h_1$ is clearly 
bounded. If $p<\infty$, since the mapping from the unit circle to itself 
induced by $B$ is $\mathcal C^\infty$ with nonvanishing derivative and 
$2$-to-$1$, we have $h_1 \in H^p(\Dd)$. 
We then write $h \in H^p (B)$. 

Since $B(0)=0$, the
composition operator $C_B$ ($C_B(f)=f\circ B$) is an isometry of $H^p$ 
into $H^p$.
Therefore, $\overline{\Cc[B]}=H^p(B)$, where 
$\Cc[B]$ stands for the polynomials in $B$, and the closure is taken 
either with respect to  $H^p$  convergence when $p<\infty$, or
to bounded pointwise convergence when $p=\infty$. 

In particular, let $g \in H^\infty(\Dd)$, then $g + g 
\circ \varphi$     and 
\[ 
 \Gamma(g)(z) := \frac{g(z) - g \circ \varphi 
(z)}{z-\varphi (z) } 
\] 
 are in $H^\infty (B)$ (the possible singularity is 
easily seen to be  removable). Furthermore, 
\[ 
g - g \circ \varphi = 2g - (g + g \circ \varphi ) 
\] 
shows that $g - g \circ \varphi  \in \Cc[H^\infty 
(B),g] \subset    \overline{\Cc[B,g]}$.

We claim that the following statement is true.
                      
                     \vspace{.2cm}

{\bf Proposition.} 
{\em Let $0<p<\infty$ and $g \in H^\infty(\Dd)$. When the clo\-sure is 
taken with respect to $H^p$, 
 $\overline{\Cc[B,g]} =$ $ \overline{\Cc[B]} +  I(g) H^p(\Dd)$ =
$H^p(B)+I(g)H^p(\Dd)$, where $I(g)$ is the inner factor 
                     in the canonical inner-outer factorization of $\Gamma 
(g)$. }

                      \vspace{.2cm}

         The proof of this result is as follows.  Let $h \in H^p(\Dd)$. 
Then 
\[ 
h (z) \Gamma(g)(z) = \frac12 (h(z) + h \circ 
\varphi(z))  \Gamma(g)(z) 
+ \frac12 \frac{h(z) - h \circ \varphi(z)}{z-\varphi (z)} 
(g(z) - g \circ \varphi(z)) . 
\] 
By applying the above 
decomposition to the terms of a sequences of 
polynomials $h_\nu$ tending to $h$ in $H^p$, 
it is then easy to see that the first term is in 
$\overline{\Cc[B]} \subset \overline{\Cc[B,g]} $, 
while $\frac{h(z) - h \circ \varphi(z)}{z-\varphi 
(z)}$ is multiplied by 
$g - g \circ \varphi$, which has been show above to be in 
$\overline{\Cc[B,g]}$.
This proves that $\Gamma(g) H^p(\Dd) \subset 
\overline{\Cc[B,g]}$, and  $\overline{\Gamma(g) H^p(\Dd)} = I(g) H^p(\Dd)$, by 
Beurling's theorem. 

Conversely, consider a monomial of the form $M= g^\al 
B^\beta$.    If $\al=0$, then $M \in \Cc[B]$. If $\al \ge 1$, 
\begin{eqnarray}
\nonumber
&&g^\al = 
\left( 
\frac12 (g - g \circ \varphi) + \frac12 (g + g \circ 
\varphi) 
\right)^\al \\
\nonumber
&&= 
\frac1{2^\al} \sum_{l=0}^\al 
\left( \begin{array}{c} \al \\ l \end{array}  \right) 
(g - g \circ \varphi)^l (g + g \circ \varphi)^{\al - l} 
\in H^\infty (B) + (g - g \circ \varphi)H^\infty (B), 
\end{eqnarray} 
 as can be seen by splitting the last sums into terms 
with $l$ 
 even or odd respectively. Now $g - g \circ \varphi = 
(z-\varphi) \Gamma(g)  \in \Gamma(g) H^\infty(\Dd)$. 

This proves that $\Cc[B,g] \subset \Cc[B] +  I(g) H^p(\Dd)$. To prove the 
corresponding inclusion for the $H^p$ closures, it will be enough to show 
that $ \overline{\Cc[B]} +  I(g) H^p(\Dd)$ is closed in $H^p$. 

Each term 
in the sum is a closed subspace. For any
$f \in H^p (\Dd)$, write
$$
f = \pi_1(f) + \pi_2(f) := \frac12 (f + f \circ \varphi) +  \frac12 (f - f 
\circ \varphi).
$$
Then the maps $\pi_1$ and $\pi_2$ are continuous from $H^p(\Dd)$ to 
itself,  $\pi_1$ is a projection onto $\overline{\Cc[B]} = H^p(B)$, and 
$H^p(B)\subset \mbox{Ker } \pi_2$. Suppose  that $f = 
\lim_n (x_n+y_n)$, with $x_n \in H^p(B)$ and $y_n \in I(g) H^p(\Dd)$. Then
$$
\pi_2(f) = \lim_{n\to\infty} \left( \pi_2(x_n) + \pi_2(y_n) \right)
= \lim_{n\to\infty} \pi_2(y_n) ,
$$
which proves the existence of the last limit. Since $\Gamma (g) \circ 
\varphi = \Gamma (g)$, we have $I (g) \circ \varphi = I (g)$, and if $y_n = 
I(g) w_n$, with $w_n \in H^p(\Dd)$, we have 
$$
\pi_2 (y_n) =  I (g) (w_n - w_n \circ \varphi) \in I(g) H^p(\Dd),
$$
therefore $\pi_2(f) \in I(g) H^p(\Dd)$ and $f \in \overline{\Cc[B]} +  I(g) 
H^p(\Dd)$, q.e.d.


\vspace{.2cm}

 
\vspace{.1cm}
{\bf Corollary.}           
{\em $\Cc[B,g]$ is dense in $H^p$ if and only if $\Gamma(g)$ is an outer 
function. }

 \vspace{.1cm}
In the case where $B(z)=z^2$ and $\Gamma(g)$ is {\it not} an outer 
function, it is easy to find an annihilator for $\overline{\Cc[B]} +  I(g) 
H^p(\Dd)$ directly. Here $\varphi (z) = -z$.
By Beurling's theorem, there exists $h \in H^q (\Dd)$ 
such that $\int_{\partial \Dd} h \overline f=0$, for any $f \in I(g) H^p(\Dd)$. 
We claim $h_1:=h-h\circ \varphi$ is the annihilator we are looking for. 
Indeed, for any $f \in I(g) H^p(\Dd)$,
$$
\int_{\partial \Dd} (h\circ \varphi) \overline f =
\int_{\partial \Dd} h \overline{ f \circ \varphi } = 0,
$$
since $f \circ \varphi \in I(g) H^p(\Dd)$. So $\int_{\partial \Dd} h_1 \overline 
f =0$.

Furthermore, if $f \in H^p(B)$, 
then $f \circ \varphi =f$ and
$$
\int_{\partial \Dd} (h-h\circ \varphi) \overline f =
\int_{\partial \Dd} h \overline f - \int_{\partial \Dd} (h\circ \varphi) \overline f =
\int_{\partial \Dd} h \overline f - \int_{\partial \Dd} h \overline{ f \circ \varphi }
=0.
$$

The case of a general Blaschke product $B$ of degree $2$ can be reduced (as far 
questions of density and the computation of $I(g)$ are concerned) to the 
case of $z^2$ by using a M\"obius map from the disk to itself taking the 
unique critical point of $B$ to $0$. We omit the details.

\vspace{.2cm}
                     
 If $B(z)=z^2$ and $g(z)=\chi (z) = 
\exp{\frac{z+1}{z-1}}$ as above, then a 
computation shows that $\Gamma(\chi)$ has no zeroes, 
and $\chi$ is continuous up 
 to the boundary and non zero for any $\zeta \in 
\partial \Dd \setminus  \{1\}$. It is then easy to see that $\Gamma(\chi)$ is 
continuous on   $\partial \Dd \setminus \{+1,-1\}$, so that any 
singular inner factor of $\Gamma (\chi )$ would have to have its singular measure concentrated 
on those two points; but  $\Gamma(\chi)$ admits a finite non-zero radial limit 
at each of those points, which precludes any singular inner factor. 
Thus, $\Cc[z^2,\chi (z)]$ is  dense in $H^2$.

 On the other hand, still for $B(z)=z^2$, if we pick 
$g(z) = z  \chi(z^2)$, then $\Gamma(g)= \chi(z^2)$, which is 
singular inner. Thus, $\Cc[B,g]$ is not dense in $H^2$. At the same time 
$z^2$ and $g(z)$ obviously
separate points of $\Dd $ (which follows also from the fact that $\Gamma 
(g)$ does not vanish inside $\Dd$).
      Moreover, by the above Proposition, for every function 
$f\in H^2$ which is orthogonal to the $z$-invariant
subspace $\chi H^2(\Dd )$ the function $h(z)=zf(z^2)$ 
is orthogonal to $\Cc[B,g]$,
so $\overline{\Cc[B,g]}$ has infinite 
codimension in $H^2$.

\vskip.5cm
{\bf Notations and statements of results}

  In this paper we study algebras with two generators, 
one of which is a finite Blaschke product
 of degree $n$:
\begin{equation}
\label{finprod}
B(z):=\prod_{j=1}^n \frac{a_j-z}{1-\bar{a}_jz},
\end{equation}
where  $a_1,...,a_n\in \Dd $. Given a point $z\in \Dd $, there are exactly 
$n$ points in the unit disk (counting multiplicities),
$z:=\varphi _0(z), \varphi _1(z),...,\varphi_{n-1}(z)$ 
such that $B(\varphi _j(z))=B(z), \ j=1,...,n-1$.
For $f\in H^\infty $ we define $\Gamma _B (f)$ by
$$
\Gamma_B 
(f)(z)=\frac{\prod_{j=1}^{n-1}(f(z)-f\circ\varphi 
_j(z))}{\prod_{j=1}^{n-1}(z-\varphi _j(z))}.
$$
In what follows we usualy omit the subscript for 
$\Gamma _B$ and write $\Gamma (f)$ instead of 
 $\Gamma_B (f)$. It will be shown in section 2 that for 
all $f\in H^\infty $, \ $\Gamma (f)$ is
 analytic and, since it is obviously bounded, is an 
$H^\infty $ function. Moreover, for all
 $0<p\leq \infty $ the mapping \ $f\longmapsto \Gamma 
(f)$ can be extended as a continuous mapping
from $H^p$ to $H^{p/(n-1)}.$  

The main result of this paper is the following
Theorem.
                
 \vspace{.3cm}

{\bf Theorem B}. 
{\em 
    Let ${\mathcal A}$ stand for the algebra 
of functions $\Cc[B,g]$ where $g\in 
H^\infty $. The closure of ${\mathcal A}$, denoted  $\overline{{\mathcal A}}$,
 is taken with respect to the 
$H^p$ norm when $p<\infty$, and with respect to bounded pointwise 
convergence when $p = \infty$.
 \begin{enumerate}

\item For every $p>0$ the closure of ${\mathcal A}$ in 
$H^p$ contains      the $z$-invariant subspace $\Gamma (g)H^p$.  
\item     $\overline{{\mathcal A}} = H^p$ if and only if 
$\Gamma (g)$ is outer. 
\item  $\overline{{\mathcal A}}$ has finite codimension 
in $H^p$
if and only if the inner part of
 $\Gamma (g)$ is a finite Blaschke product.
\end{enumerate} }
                      
\vspace{.2cm}

As a corollary to the proof of Theorem B we obtain 
the following  generalization of Theorem A 
(though $p_1$ is not a Blaschke product,  $\Gamma 
_{p_1}(p_2)$ below can be defined 
in a similar way. The details may be found in 
section 4). 

\vspace{.3cm}

{\bf Theorem C}. {\em Let $p_1$ and $p_2$ be two 
polynomials. If $\Gamma_{p_1}(p_2)$
is not identically equal to zero, then the algebra 
$\overline{\Cc[p_1,p_2]}$ has  finite
codimension in $H^2$. Furthermore, $\Cc[p_1,p_2]$ is dense in $H^2$ if and
only if $\Gamma_{p_1}(p_2)$ does not vanish in the open unit disk.}

 \vspace{.2cm}
                     
In the case when $\Gamma (g)\equiv 0$, Theorem B (as 
well as Theorem C)  does not give much information about the size of the 
closure of the algebra. It is shown in section \ref{sec vanish} below that 
the identical vanishing of $\Gamma (g)$
 implies a strong functional relation between $B$ and 
$g$. The following  Theorem gives  a necessary and sufficient condition 
for the identical  vanishing of $\Gamma (g) $.
              
\vspace{.3cm}

{\bf Theorem D}. {\em The function $\Gamma (g)$ vanishes identically in 
the unit disk if and only if there exist a function $\varphi$ from the disk to
itself and a discrete subset $A$ of the disk such that for any $z \in \Dd
\setminus A$, $\# \varphi^{-1} (z) > 1$, such that $B = B_1 \circ \varphi$ 
and $g =g_1 \circ \varphi$,  where $B_1$ and $g_1$ are functions defined on the disk.}
\vspace{.2cm}

 A direct corollary to this result is 

 \vspace{.2cm}

{\bf Theorem $\mbox{D}^\prime $}. {\em The 
function $\Gamma (g)$ vanishes identically if and only
if the closure of $\Cc[B,g]$ does not contain a 
non-trivial $z$-invariant subspace.}

\vspace{.2cm} 

As another corollary to  Theorem D we obtain the 
following result.

 \vspace{.2cm}

{\bf Corollary}. {\em If the order of the 
Blaschke product $B$ is a prime number,
then, whenever the bounded analytic function $g$  
 is not in $H^\infty [B]$, \ $\Gamma (g)$ does not vanish identicaly.}

\vspace{.2cm}

Finally, though the disk algebra $A(\bar{\Dd})$ 
is not our prime focus, our methods
yield similar disk algebra results. It is easily 
seen that if $g$  is a disk algebra function, $\Gamma (g)$ is also in 
the disk algebra.    The following Theorem holds.

\vspace{.3cm}

{\bf Theorem E}. {\em 
\begin{enumerate}
\item Let $g\in A(\bar{\Dd})$. Then the closure 
of $\Cc[B,g]$  in $A(\bar{\Dd})$ contains 
$\Gamma(g)A(\bar{\Dd})$.
 \item $\Cc[B,g]$ is dense in $A(\bar{\Dd})$ if 
and only if $\Gamma (g)$ has no zeros in the closed unit disk.
\item $\overline{\Cc[B,g]}$ has finite 
codimension in $A(\bar{\Dd})$ if and only if 
$\Gamma (g)$ has no singular inner factor and 
only a finite number of zeros in the closed unit disk.
\end{enumerate} }
 
\vspace{.2cm}

The structure of this paper is as follows. 
Section 2 is devoted to some auxiliary results related to level sets of the 
Blaschke product $B$. Here we introduce the function $\Gamma (g)$ and 
prove its analyticity. In section 3 we prove the ``necessary" part of Theorem 1. In 
section 4 we show how Theorem C  follows from Theorem B.  In  this section  we 
also give an estimate of the number of generators of a closed submodule  of $H^2$ 
over a polynomial algebra expressed in terms of the function $\Gamma $. Section 5 
deals with the case when $\Gamma (g) \equiv 0$.  Theorems D and $\mbox{D}^\prime $ are proved in 
this section. In section 6 we discuss $z$-invariant subspaces contained in 
the algebra's closure.   Section 7 is devoted to applications to the 
disk algebra case.
Here we also compare our results with  previous 
ones obtained by J. Wermer and R. G. Blumenthal
without the assumption that one of the generators be a finite Blaschke product.
Finally, in section 8 we conclude  the proofs of Theorems B and E.

{\bf Acknowledgements}. The authors would like 
to thank John Wermer for stimulating discussions and the references 
\cite{JJo}, \cite{JWe}, and Emma Previato for pointing out and explaining 
\cite{JLB}.
We thank the referee for his careful comments and corrections.

\section{Auxiliary results}
\label{sec sep}
 \subsection{Level sets of $B$.} 
\label{level sets}

For any $z \in \overline{\Dd}$, the set $B^{-1}(B(z))$ can be written as 
$\{z=: \varphi_0 (z), \varphi_1 (z), \ldots , 
\varphi_{n-1} (z) \}$, 
which are the $n$ solutions (possibly equal) of a 
polynomial equation of 
degree $n$; 
the maps $\varphi_j$ may be chosen holomorphic in $z$ near any point $z^0$ such that 
$\varphi_j (z^0)$ is a simple zero of $B - B(z^0)$ (equivalently, such that $ B' 
(z^0) \neq 0$), but cannot be chosen to be
globally continuous. 

Notice in particular that $B'$ cannot vanish on the 
boundary  $\partial \Dd$ of the disc, and that, since they all 
map the circle into itself and never coincide, the maps $(\varphi_j, 1 \le 
j \le n-1)$ can be ordered by circular order, 
and thus are all individually
well-defined in a neighborhood of the unit circle. Explicitly, we can 
write $\varphi_j (\theta) = \exp(i 
\alpha_j(\theta))$, where $\theta < \alpha_j(\theta) < \theta + 2 \pi$. A 
priori, those choices only hold locally.
Then 
label the $n-1$ distinct points of the unit circle $\{ \varphi_1 (1), \ldots , 
\varphi_{n-1} (1) \}$ by demanding $\alpha_j(0) < 
\alpha_{j+1}(0)$, $1 \le j \le n-2$. Because $B'(\eit)\neq 0$ for 
any $\theta$, we never have $\lim_{\theta\to\theta_0}\alpha_j(\theta) = 
\lim_{\theta\to\theta_0}\alpha_{j+1}(\theta)$, so the set
$$
U_j := \{ \theta \in [0,2\pi) : \alpha_j, \alpha_{j+1} \in {\mathcal {C}}^0 
[0,\theta] \mbox { and }
\forall \theta' \in [0,\theta], \alpha_j(\theta) < 
\alpha_{j+1}(\theta) \}
$$
is both open and closed, so is equal to the whole interval $[0,2\pi)$ and 
passing to the limit and using $B'(\eit)\neq 0$ again, we will have 
$$
2\pi < \alpha_{j}(2\pi) < 
\alpha_{j+1}(2\pi) <  4 \pi, 1 \le j \le n-2,
$$
and therefore $\alpha_{j}(2\pi) = \alpha_{j}(0) + 2 \pi$, $0 \le j \le 
n-1$, and we have liftings of continuous maps.

Furthermore, the quantities $|\varphi_j
 (e^{i\theta})- \varphi_k (e^{i\theta})|$, $0\le \theta \le 2 \pi$, $0\le j \neq k \le n-1$, are 
uniformly bounded away from $0$. 

Near a multiple point $z^0$, 
so that $B-B(z_0)$ has a zero of order $N$, the 
$\varphi_j$   can be expressed as series in $(z-z^0)^{1/N}$ (so 
they are actually continuous at the point $z^0$, if not in any neighborhood of it) 
(cf.  for example \cite{AIM}).   

                     Symmetric polynomial (resp. rational) functions of 
$\varphi_j (z), 0 \le j \le n-1$ 
                     are well-defined and holomorphic (resp. meromorphic) 
on the whole disc.   
                     For a multiple point $z^0$, if $A \subset \{0, 
\ldots, n-1\}$ is such that 
                     $\varphi_j (z^0)= \varphi_k (z^0)$ for any $j$, $k 
\in A$, and 
                     $\varphi_j (z^0)\neq \varphi_k (z^0)$ for any $j\in 
A$, $k \notin A$, 
                     then the symmetric polynomial functions of 
$(\varphi_j (z), j \in A)$ 
                     are well-defined and holomorphic in some neighborhood 
of $z^0$ 
                     (as can be seen by factoring the equation). 

                     In the special case that is of interest to us, since 
                     $\varphi_0 (z) =z$ is always globally well-defined, 
the same is true 
                     of the symmetric functions of $(\varphi_j (z), 1 \le 
j \le n-1)$. More precisely, 
 we have the following easy lemma. 
From now on, we shall always assume
that $B$ vanishes at the origin, so in (\ref{finprod}) we have
 $a_1=0$. This is no loss of generality as 
it can be  obtained by composing our functions by a M\"obius map 
of the disk. 

\begin{lemma} 
\label{sym_func} 
Any symmetric polynomial (resp. rational function, resp. holomorphic 
function in the disk)  in 
$\varphi_j, 0 \le j \le n-1$, can be written as a polynomial (resp.
rational function) in $B(z)$. 

Any symmetric polynomial (resp. rational function)  in 
$\varphi_j, 1 \le j \le n-1$, can be written as a polynomial (resp.
rational function, resp. holomorphic 
function in the disk) in $z$ and $B(z)$. 
\end{lemma} 

\begin{proof} 
Denote the elementary symmetric functions in $n$ 
(resp $n-1$) indeterminates by 
\begin{eqnarray*} 
\sigma_k (X_0, \ldots X_{n-1}) &:=& 
\sum_{ A \subset \{0, \ldots, n-1\}, \# A=k} \, \prod_{j 
\in A} X_j , 
\\ 
{\rm resp.\ } 
\sigma'_k (X_1, \ldots X_{n-1}) &:=& 
\sum_{ A \subset \{1, \ldots, n-1\}, \# A=k} \, \prod_{j 
\in A} X_j . 
\end{eqnarray*} 
It is clear that the $\sigma_k (\varphi_0 (z), 
\ldots,\varphi_{n-1}(z))$ depend only on $B(z)$.  The polynomial equation 
which admits  $\{ \varphi_0 (z), \ldots , \varphi_{n-1} (z) \}$ as 
solutions  is 
 \[ 
\prod_{j=1}^n (a_j-X) = B(z) \prod_{j=1}^n (1-\bar a_j X). 
\] 
Since $a_1=0$, the coefficient of $X^n$ is exactly 
$(-1)^n$, and all the coefficients of $X^m$, $0 \le m \le n$, 
are affine functions of $B(z)$, which proves the first 
assertion. 

To prove the second assertion of the lemma, it will be 
enough to show that the $\sigma'_k$ are polynomials in $X_0$ and
$\sigma_k$. But 
\[ 
\sigma_k (X_0, \ldots X_{n-1}) = 
X_0 \sigma'_{k-1} (X_1, \ldots X_{n-1}) + 
\sigma'_k (X_1, \ldots X_{n-1}), 
\] 
with $\sigma_0 = \sigma'_0 =1$, so the relation 
$\sigma'_k = \sigma_k - X_0 \sigma'_{k-1}$ allows us to conclude by
induction. 
\end{proof} 

\vskip.5cm 

\subsection{Auxiliary functions.} 

For any $f \in \calo (\Dd)$, define 
\[ 
\cald (f)(z) := \prod_{j=1}^{n-1} 
\left(f(z)-f\circ\varphi_j (z)\right). 
\] 
This is symmetric in $(\varphi_j (z), 1 \le j \le 
n-1)$, and, therefore, by Lemma 1 is globally holomorphic in $z$. In particular, 
$\cald(z) := \prod_{j=1}^{n-1} \left(z-\varphi_j 
(z)\right)$ is a function which 
vanishes at $z^0$ if and only if $B'(z^0)=0$ (it can be seen to be a discriminant). 

Define $\Gamma(f)(z) := \cald(f)(z)/ \cald(z)$. This 
is clearly meromorphic in the disc, but actually each 
singularity is removable. 
Indeed, near any point $z^0$ where for some $j$, 
$\varphi_j(z^0)=z^0$, 
 we have 
\[ 
f(z)= \sum_{m \ge 0} \alpha_m (z-z^0)^m , \quad 
f(\varphi_j(z))= \sum_{m \ge 0} \alpha_m 
(\varphi_j(z)-z^0)^m , 
\] 
where the second series converges in a neighborhood of $z^0$ 
since $|\varphi_j(z)-z^0| $ 
\newline
$= |\varphi_j(z)-\varphi_j(z^0)|  \le C |z-z^0|^{1/N}$. Therefore, 
\begin{eqnarray}
\nonumber
& f(z) - f(\varphi_j(z)) =
\\ 
\nonumber
& (z-\varphi_j(z)) \sum_{m \ge 1} 
\alpha_m \left( (z-z^0)^{m-1} + (z-z^0)^{m-2} 
(\varphi_j(z)-z^0) + \cdots + (\varphi_j(z)-z^0)^{m-1} \right) , 
\end{eqnarray}
and the infinite series converges in a neighborhood 
of $z^0$ because its general term is bounded by 
$C m |\alpha_m | |z-z^0|^{(m-1)/N}$. Let $A := \{ j \ge 1 : \varphi_j 
(z^0)=z^0\}$, 
 then 
\newline
$\prod_{j \notin A} \left(f(z)-f\circ\varphi_j (z)\right)/ 
\left(z-\varphi_j (z)\right)$ 
is clearly holomorphic in a neighborhood of 
                     $z^0$, and 
                     \begin{eqnarray}
\nonumber
&& \prod_{j \in A} 
\frac{\left(f(z)-f\circ\varphi_j (z)\right)} 
{\left(z-\varphi_j (z)\right)} 
 \\ 
\nonumber
&& = \prod_{j \in A} \sum_{m \ge 1} 
\alpha_m \left( (z-z^0)^{m-1} + (z-z^0)^{m-2} 
(\varphi_j(z)-z^0) + \cdots + (\varphi_j(z)-z^0)^{m-1} \right). 
\end{eqnarray}
The coefficients of the product series are easily seen 
to be symmetric in $(\varphi_j (z), j \in A)$, therefore 
well-defined and  holomorphic in a neighborhood of $z^0$, so the 
product on the left hand side of the equation, and therefore the 
function $\Gamma(f)$, has a removable singularity at $z^0$, as required. 

The zeros of the function $\Gamma(f)$ have a simple 
interpretation in terms of self-inter\-sec\-tion or cusps of the 
analytic graph   $\{(B(\zeta),f(\zeta)), \zeta \in \Dd\} 
\subset \Cc^2$. 

\begin{lemma} 
\label{zeroGamma} 
Let $z \in \Dd$, then $\Gamma (f)(z)=0$ 
if and only if   one of the following occurs: 

(i) There exists $z'\neq z$ such that $B(z')=B(z)$ 
and  $f(z')=f(z)$; 

(ii) $B'(z)=f'(z)=0$. 
\end{lemma} 

\begin{proof} 
The product defining $\Gamma(f)$ vanishes if and only 
if one   its factors does, i.e. there exists $j$, $1 \le j \le 
n-1$,    such that 
\[ 
\frac{f(z)-f\circ\varphi_j (z)} 
{\left(z-\varphi_j (z)\right)} =0. 
\] 
If $z \neq \varphi_j (z)$, this simply means that 
$f(z)=f(\varphi_j (z))$, and we have case (i). If $z =\varphi_j (z)$, $B$ is 
not locally   invertible at $z$, therefore $B'(z)=0$, and 
the vanishing of the factor means that 
\[ 
0= \lim_{z'\to z} \frac{f(z')-f\circ\varphi_j (z')} 
{\left(z'-\varphi_j (z')\right)} = 
\lim_{z'\to z} \int_0^1 f'(tz'+(1-t)\varphi_j (z')) \, dt , 
\] 
which happens exactly when $f'(z)=0$. 
\end{proof} 

\section{Necessary conditions for density or finite codimension}
\label{sec nec cond}
 We shall now see how the function $\Gamma(g)$ defined 
above
relates to the properties of the algebra $\Cc[B,g]$. The main result of 
this section is Theorem \ref{NCGamma}.

\begin{lemma} 
 \label{contGamma} 
For $0<p\le \infty$, the mapping defined by $f \mapsto \Gamma(f)$ is 
continuous from $H^p$ to $H^{p/(n-1)}$. 
\end{lemma} 

\begin{proof} 
We have seen that $\Gamma(f)$ is holomorphic on 
the whole  disc. 
By the remarks at the beginning of \S \ref{level sets}, in a neighborhood of the unit circle 
the maps $\varphi_j$ are all globally defined, holomorphic, 
with a non-vanishing derivative, it is clear that each function 
 $f \circ \varphi_j$ is an $H^p$ function : holomorphic on a neighborhood
of the unit circle in $\Dd$, bounded for $p = \infty$, and
for $p<\infty$,  $f \circ \varphi_j (r \zeta)$ tends to $f \circ 
\varphi_j (\zeta)$  in the $L^p$ norm on the unit circle. The $H^p$-norm 
of $f \circ \varphi_j$ is
bounded by  $C \|f\|_{H^p} = C \|f\|_{L^p(\partial \Dd)}$. 

Since $|z-\varphi_j (z)|$ is bounded from below in a 
neighborhood of the unit circle, a similar $H^p$ bound holds for $(f-f 
\circ \varphi_j)/(z-\varphi_j (z))$.  Now the map $\Gamma$, seen as a map from $H^p$ to 
$L^{p/(n-1)}(\partial \Dd)$, is a composition of linear and $(n-1)$-linear 
maps, so to check   that it is continuous it is enough to see that 
\[ 
\|\Gamma(f)\|_{L^{p/(n-1)}(\partial \Dd)} 
\le C \|f\|_p^{n-1}, 
\] 
 which follows from the above estimates and repeated 
use of H\"older's inequality. 

For $p=\infty$, the required estimate follows 
immediately from $|z-\varphi_j (z)|$ being bounded from below. 
\end{proof} 

 \begin{proposition} 
\label{Aincl} 
Let $\cala := \Cc [B,g]$, and $\overline {\cala}$ be 
the closure of $\cala$ with respect to the $H^p$ norm, for some $p 
\in (0,\infty)$. Let $\Gamma(g)= B_g S_g F_g$ be the canonical 
factorization of $\Gamma(g)$ into, respectively, Blaschke product, 
singular inner factor, and outer factor. Then for any $f \in \overline 
{\cala}$, $\Gamma(f) \in B_g S_g H^{p/(n-1)}$. 
\end{proposition} 

\begin{proof} 
Suppose first that $f \in \cala$. Then, using the 
fact that $B \circ \varphi_j = B$ for any $j$, 
\begin{eqnarray*} 
f&=& \sum_{\al, \beta \ge 0} a_{\al \beta} g^\al 
B^\beta \mbox{ \ implies} \\ 
& &\Gamma(f) = \Gamma(g) \, 
\prod_{j=1}^{n-1} \left( 
\sum_{\al \ge 1, \beta \ge 0} a_{\al\beta} B^\beta 
\left[ \sum_{l=0}^{\al-1} g^{\al-1-l} (g \circ 
\varphi_j)^l \right] 
                     \right), 
\end{eqnarray*} 
                     so $\Gamma(f) \in \Gamma(g) \Cc[H^\infty(B),g]$, 
because the last product involves a finite set of powers of $B$ and 
polynomials which are symmetric in the $\{g \circ \varphi_j, 1 \le 
n-1\}$, and therefore  can be expressed as polynomials of $g$ and bounded 
holomorphic functions  symmetric in $\{g \circ \varphi_j, 0 \le n-1\}$, 
which are then necessarily bounded holomorphic functions of $B$. 

                   Suppose $\cala \ni f_\nu \to f$ in $H^p$. By the above 
we have 
                     $\Gamma (f_\nu) = \Gamma(g) q_\nu$, $q_\nu \in 
\Cc[H^\infty(B),g]$. 
                     Now in the canonical factorization, the harmonic 
function 
                     $\log |F_g|$ is the Poisson integral of its boundary 
values, 
                     while $\log |S_g|$ is the Poisson integral of a 
singular measure 
                     $- \sigma_g$, where $\sigma_g$ is a positive measure 
(carried on a 
                     Borel set of Lebesgue measure $0$). Factoring both 
$\Gamma (f_\nu)$ 
                     and $q_\nu$, and using the essential uniqueness of 
the factorization, 
                     we see that 
                     \[ 
                     \sigma_{f_\nu} \ge  \sigma_g 
                     \] 
                     for all indices $\nu$ (their difference is a positive 
measure). 
                     Since $\Gamma (f_\nu)$ must converge to $\Gamma (f)$ 
in $H^{p/(n-1)}$, 
                     $ \sigma_{f_\nu}$ converges weakly to $ \sigma_{f}$, 
so 
                     $ \sigma_{f} \ge \sigma_g$, which implies that $S_g$ 
                     divides $\Gamma(f)$. 

                     On the other hand, $H^{p/(n-1)}$ convergence implies 
convergence of the 
                     function and all its derivatives uniformly on any 
compactum, and 
                     therefore since each $\Gamma (f_\nu)$ has zeroes at 
the points 
                     of $B_g^{-1}(0)$ with at least the same multiplicity, 
the same holds true 
                     for $\Gamma (f)$. Since dividing by $B_g S_g$ only 
affects the boundary 
                     values of $\log | \Gamma(f) |$ on a set of measure 
$0$, the remaining 
                     function is still in $H^{p/(n-1)}$, q.e.d. 
                     \end{proof} 

                     \begin{theorem} 
                     \label{NCGamma} 
                     \begin{enumerate} 
                     \item 
                     If $\cala$ is dense in $H^p$, then $\Gamma(g)$ is an 
outer function. 
                     \item 
                     If 
                     $\overline {\cala}$ has finite codimension in $H^p$, 
then the inner part of 
                     $\Gamma(g)$ is a finite Blaschke product (and in 
particular there is no 
                     singular inner factor). 
                     \end{enumerate} 
                     \end{theorem} 

                     \begin{proof} 
                     If $z \in \overline {\cala}$, then by Proposition 
\ref {Aincl}, $1= \Gamma(z) \in B_gS_g H^{p/(n-1)}$, 
                     which implies that $\Gamma(g)$ has no inner factor by 
Beurling's theorem. 

                     If $\overline {\cala}$ has finite codimension $m$, 
then it must contain 
                     a non-zero polynomial $P$ of degree $m$. Then 
$\Gamma(P)$ is analytic 
                     across the unit circle, so cannot admit a singular 
inner factor. On the 
                     other hand, any zero $a$ of $\Gamma(g)$ in $\Dd$ 
provides a continuous 
                     functional on $H^p$ which vanishes on $\overline 
{\cala}$ : in the case (i) of 
                     Lemma \ref{zeroGamma}, the functional is $f \mapsto 
f(a)-f(a')$, in the case (ii), 
                     $f \mapsto f'(a)$. It is easy to see that if the set 
of zeros of 
                     $\Gamma(g)$ is infinite, this family of functionals 
is not of finite 
                     dimension, which contradicts the hypothesis. 
                     \end{proof} 

\section{Main results: sufficient conditions for density, generating sets 
for submodules}
\label{sec main}

                     \begin{theorem} 
                     \label{main claim} 
                     For any $p>0$, 
                     $\overline{\cala}$ contains 
                     $\Gamma (g) H^p$. 
                     \end{theorem}

                     Let us defer the proof of Theorem \ref{main claim}
to Section \ref{sec proof} in favor of a few observations. 

             First, by Beurling's theorem again, if $\Gamma(g)$ is outer, 
then the closure of $\Gamma (g) H^p$ is $H^p$ itself; 
and if the inner factor of $\Gamma(g)$ reduces to a finite 
Blaschke product $B_g$, then the closure of $\Gamma (g) H^p$ is 
$B_g H^p$. Since a function $f \in H^p$ is divisible by $B_g$ if and only if 
it vanishes on the zeroes of $B_g$ to an 
order greater or equal than that of $B_g$, $B_g H^p$ is then a 
subspace of finite codimension. 

                     So Theorem \ref{main claim} has as an immediate 
consequence the converse statements to those in Corollary \ref{NCGamma},
and together they make up Theorem B from the introduction. 

                     \begin{corollary} 
                     \label{NaSCGamma} 
                     \begin{enumerate} 
                     \item 
                     $\cala$ is dense in $H^p$ if and only if $\Gamma(g)$ 
is an outer function. 
                     \item 
                     $\overline {\cala}$ has finite codimension in $H^p$ 
if and only if the inner part of 
                     $\Gamma(g)$ is a finite Blaschke product (and in 
particular there is no 
                     singular inner factor). 
                     \end{enumerate} 
                     \end{corollary}

                     An application of Theorem \ref{main claim} to the 
polynomial case gives the following extension of the results proved in \cite{JLS}. 

                     Let $p(z)$ and $q(z)$ be two polynomials. An argument 
similar to the one above alows to define ${\mathcal D} _p(q)$ and $\Gamma _p(q)$ by 
$$ 
\begin{array}{l} 
{\mathcal D}_p(q)(z):=\prod_{j=1}^{n -1} 
(q(z)-q(\varphi_j(z))  
\\ 
\Gamma _p (q)(z):={\mathcal D}_p(q)(z)/{\mathcal D}_p(z) \\ 
\end{array} 
$$ 
                     where again $\{ z=\varphi _0(z), \varphi 
_1(z),...,\varphi _{n-1}(z) \}$ are $n$ roots 
                     of the polynomial equation (in $t$) $p(t)=p(z)$, and 
$n=\mbox{deg}\ p$. The same argument shows that 
                     ${\mathcal D}_p(q)$ and $\Gamma _p(q)$ are polynomials. 

The following result has Theorem C as a special case.

\begin{theorem} 
\label{Polynom} 
Let $p_1,...,p_m\in \Cc[z]$. Suppose that for some 
$1\leq j\leq n$, at
 least one of $\Gamma_{p_j}(p_k), \ k=1,2,...,m, \
k\neq j$, is not identically zero. 
Then $\overline{\mathcal A}=\overline{\Cc[p_1,...,p_m]}$ 
has finite codimension in $H^2$, and every closed 
${\mathcal A} $-submodule of $H^2$ has at most
 $\min_{k\neq j} \{ {\mathcal N} (\Gamma_{p_j}(p_k) ) \} $ generators, 
where for a polynomial $R$, ${\mathcal N} (R)$ is the
 number of zeros of $R$ in the unit disk. 
\end{theorem}

                     \begin{proof} 
                     Let $M$ be a positive number large enough so that the 
interior $\Omega  $ of some component of the lemniscate $\{ z:|p_j(z)|= M \} $ 
contains the unit disk, and $p'(z)\neq 0$ whenever
$z \in \partial \Omega$. Denote by $\psi $ the  Riemann mapping from $\Omega $ onto 
the  unit disk , which maps $p_j(0)$ to the origin. 
Then $\psi $ is analytic in $\bar{\Omega }$ and the 
Hardy space in $\Omega $, $H^2(\Omega )$ (we equip $\partial \Omega $ with harmonic measure), 
is isometric to $H^2$ in the unit disk: an isometry is being given by 
$$ 
F: f\longmapsto f\circ \psi . 
$$ 
Note that $F(p_j)$ is a finite Blaschke product and 
$F(p_l), \ l=1,...,m$, are analytic in a neighborhood of 
the unit disk. This implies that $\Gamma 
_{F(p_j)}(F(p_k))=F(\Gamma _{p_j}(p_k))$ have trivial singular 
parts in their inner-outer factorization for all $k=1,\dots,m$ and their 
inner parts are finite 
 Blaschke products. By Theorem \ref{main claim}, $F(\overline{\mathcal A})$
contains $z$-invariant subspaces generated 
 by these Blaschke products. In particular this implies that for every 
fixed $k$, \ ${\mathcal A}$ contains all  polynomials vanishing with greater or equal 
multiplicity
 at all zeros of $\Gamma _{p_j}(p_k)$ in $\Omega $. Since the $H^2(\Omega 
)$ convergence 
 implies $H^2$-convergence in the unit disk, Lemma 3 of \cite{JLS} implies 
that the
closure of 
 ${\mathcal A}$ in $H^2$ contains the $z$-invariant subspace determined by the
zeros of 
$\Gamma _{p_j}(p_k)$ which lie in the interior of the unit disk,
 and, therefore, has finite codimension. 
    Let $B_k$ be the inner part of $\Gamma _{p_j}(p_k)$. 
Thus, since every  ${\mathcal A}$-submodule $S$ of $H^2$ is a 
$\Cc[B_k]$-submodule, the
number of generators of $S$ 
as ${\mathcal A}$-module does not exceed the number of generators of $S$ as
$\Cc[B_k]$-module. By the Wold 
decomposition theorem this last number does not exceed the order of $B_k$ 
(see
\cite{JLS} for the details). 
                     \end{proof} 

                     \begin{corollary}
                      If $\Gamma _{p_j}(p_k), \ k=1,...,m, \ k\neq j$ have 
no common zeros in the unit disk, then
                      ${\mathcal A}=\Cc[p_1,...,p_m]$ is dense in $H^2$.
                      \end{corollary}

\begin{proof} Let $a_1,...,a_N$ be  common zeros in 
$\Cc$ of the   polynomials  $\Gamma_{p_j}(p_k)$. There are 
polynomials  $ q_k, \ k=1,...,m, \ k\neq j$ such that 
$$
\sum_k \Gamma _{p_j}(p_k)(z)q_k(z) =R(z)=\prod_{l=1}^N (z-a_l).
$$
                      By Theorem \ref {main claim} $\overline{R(z)H^2}$ is 
in $\overline{{\mathcal A}}$. Since 
                      $|a_l|\geq 1$, \ $R(z)$ is outer, and 
$\overline{R(z)H^2}=H^2$. 
\end{proof}

Note that in general the estimate given by Theorem 
\ref{Polynom} is sharp. As an example 
consider the algebra ${\mathcal A}=\Cc [p_1,p_2]$ where 
$p_1(z)=z(z-a), \ p_2(z)=z^2(z-a)$ 
for some $|a|<1$.  It is easily seen that $\Gamma 
_{p_1}(p_2)(z)=p_1(z)$. 
Thus, by Theorem \ref{Polynom} every closed ${\mathcal 
A}$-submodule of $H^2$ has at most two generators. 
Let $M$ be the subspace of $H^2$ consisting of all 
functions vanishing at some point $b$ in the unit disk,
different from $a$ and $0$. 
Clearly $M$ is a closed ${\mathcal A}$-submodule which is 
not singly generated. Indeed, for every function $g\in M$
such that $g(0)\neq 0$ and $g(a)\neq 0$,  there exists a function 
$h\in M$ such that  $g(0)h(a)\neq h(0)g(a)$, and, therefore, $g$ does not 
generate $M$ as ${\mathcal A}$-module. If $g(0)=0$ or $g(a)=0$, then since for 
any $f\in \cala$, $f(0)=f(a)=0$, the $\cala$-module generated by $g$ 
cannot contain $z-b$, for instance.

\section{Case $\Gamma (g)\equiv 0$.}

\label{sec vanish} 

At the other end of the spectrum from the density 
cases, there are cases when $\overline{\cala}$ cannot 
contain any non-trivial   invariant subspace --- the simplest case would be 
when $B$ and $g$ are both even functions, for instance: then any function in 
$\overline{\cala}$ is also even, and multiplying by $z$ takes us out of that subset.  In 
that case, naturally, $\Gamma(g)$ vanishes identically, so Theorem \ref{main claim} 
becomes vacuous. One would then like to know when $\Gamma(g) \equiv 0$. 

A typical example of a situation when $\Gamma(g) 
\equiv 0$ is when there 
exists a function $u $ such that 
$B=B_1 \circ u$ and $g=g_1 \circ u$, and such that 
there exists a non-empty open set $\Omega \subset \Dd$ such that for 
any $z \in \Omega$,  there exists $z'\neq z$, $z'\in \Dd$, such that 
$u(z')=u(z)$. This implies 
 that $B(z')=B(z)$ and $g(z')=g(z)$, therefore 
$\Gamma(g)$ vanishes on $\Omega$, 
and thus on $\Dd$. 

Theorem D states that this is how the identical 
vanishing of $\Gamma(g)$ always arises; it follows 
immediately from the following theorem. 

Let us introduce a notation for 
the (finite) set of $z$ such that 
$\{ \varphi_j (z), 0\le j \le n-1\}$ is made up of 
fewer than $n$ points: 
\begin{equation} 
\label{defso}
 \cals_0 := \{ z: B'(\varphi_j (z))=0 {\rm\ for\ some\ } 
0\le j \le n-1 \}. 
\end{equation} 

 \begin{theorem} 
\label{vanishGamma} 
The function $\Gamma(g)$ vanishes identically on the 
unit disk if and only if 
there exists an open set $\Omega \subset \Dd 
\setminus \cals_0$ such that $\Dd \setminus \Omega$ 
                     is discrete, and an integer $m \ge 2$ 
                     such that for any $z \in \Omega$, the 
                     set $B^{-1}(B(z))$ is partitioned into subsets $A_p$ 
of cardinality $m$ such 
                     that $g$ is constant on each $A_q$. 
                     \end{theorem} 

Before passing to the proof of this theorem, notice that this means in particular that $m$ must 
divide $n$, therefore we obtain the following

\begin{corollary}
If $n$ is prime, $\Gamma(g)\equiv 0$ implies that 
$g$ must in fact be a function of $B$, and $\overline{\cala} = 
\overline{\Cc[B]}$.
\end{corollary} 
                                    
Another immediate corollary to Theorem \ref 
{vanishGamma} is the following result.
  
\begin{corollary}
\label{vanish ident} 
The function $\Gamma(g)$ vanishes identically in the 
unit disk if and only if 
the only $z$-invariant subspace of $H^p$ contained in 
$\overline{\cala}$ is the trivial subspace $\{0\}$. 
\end{corollary} 

\begin{proof} 
If $\Gamma(g) \not \equiv 0$, a non-trivial 
$z$-invariant subspace contained in $\overline{\cala}$ 
is given by Theorem \ref{main claim}.   

If $\Gamma(g) \equiv 0$ and $M$ is a 
nontrivial $z$-invariant subspace contained in 
$\overline{\cala}$, choose $f \in M \setminus \{0\}$. Choose a point $a 
\in \Omega \setminus \{0\}$ such that $f(a) \neq 0$; then for any $b\neq a$ such 
that $B(b)=B(a)$ and $g(b)=g(a)$
(the above Theorem guarantees the existence of such points), 
we must have $f(b)=f(a)$, but $bf(b)\neq af(a)$, 
which shows that $zf \notin \overline{\cala}$, 
thus $zf \notin M$ and $M$ is not $z$-invariant. 
\end{proof} 

                     \begin{proof}[Proof of Theorem \ref{vanishGamma}] 
                     Define an integer-valued function in the disk by 
                     \[ 
                     C(z):= \# \{ w \in B^{-1}(B(z))  : g(w)=g(z) \}. 
                     \] 
                     The hypothesis says that 
                     $C(z) \ge 2$ for any $z \in \Omega_0 := \Dd \setminus 
\cals_0$. 

                     We claim that $C$ is upper-semicontinuous on 
                     $\Omega_0$. 

                     Indeed, let $z_0 \in \Omega_0$, and let $\{z_\nu\}$ a 
sequence tending 
                     to $z_0$ such that $\lim_{\nu\to\infty} C(z_\nu)= 
\limsup_{z \to z_0} 
                     C(z)$. We may assume that $\{z_\nu\}$ 
                     is contained in a neighborhood of $z_0$ 
                     on which $\varphi_j$, $1\le j \le n-1$ are all 
well-defined and distinct. 
                     For each $\nu \in \Zz_+$, there exists a subset 
$A_\nu \subset \{0,\dots 
                     n-1\}$ such that $g \circ \varphi_j(z_\nu) = 
g(z_\nu)$ if and only if $j \in 
                     A_\nu$, and $\# A_\nu = C(z_\nu)$. For $\nu$ large 
enough, this is 
                     constant and equal to $\limsup_{z \to z_0} C(z)$. 
Since there is only a 
                     finite number of subsets of $\{0,\dots n-1\}$, we may 
also assume, by 
                     passing to a subsequence again denoted by 
$\{z_\nu\}$, that $A_\nu=A$ is 
                     a constant subset. Therefore for any $j \in A$, 
                     \[ 
                     g \circ \varphi_j (z_0) = 
                     \lim_{\nu\to\infty} g \circ \varphi_j (z_\nu) = 
                     \lim_{\nu\to\infty} g (z_\nu) = g(z_0). 
                     \] 
                     So $C(z_0) \ge \limsup_{z \to z_0} C(z)$, as claimed. 

                     Let $m:= \min_{z \in \Omega_0} C(z) \ge 2$. This 
makes sense since $C$ is 
                     integer-valued, and since it is u.s.c. the set 
\[ 
\Omega_1 := \{ z \in \Omega_0 : C(z)=m \} 
\] 
is open. 

 We claim that $E:=\Omega_0 \setminus \Omega_1$ (and 
therefore $\Dd \setminus 
\Omega_1$) only has isolated points. 

Indeed, let $E'$ be the set of non-isolated points 
within $E$. Since $E$ 
is closed in $\Omega_0$, so is $E'$. On the other 
hand, suppose $z_0 \in E'$. Restrict attention to a neighborhood
of $z_0$ where a continuous choice of $\varphi_j$ can be made.
Then, reasoning as in the proof of the previous 
claim, we can find a 
sequence $\{z_\nu\} \subset E$ tending to $z_0$ and a 
set $A\subset \{0,\dots n-1\}$ such that $g \circ 
\varphi_j(z_\nu) = g(z_\nu)$ 
 if and only if $j \in A$. Then, for any $j \in A 
\setminus \{0\}$,   the analytic function $g-g \circ \varphi_j$ vanishes 
on a subset which   has a limit point, and therefore is identically $0$ in 
a small disk $z_0$. This proves that any point in this disk is in 
$E$, and therefore in $E'$. So $E'$ is also an open set. Since 
$\Omega_0$ is connected, if 
$E' \neq \emptyset$, then $E' = \Omega_0$, which is 
absurd. Therefore  $E'=\emptyset$, q.e.d. 

                     Now let 
\[ 
\Omega_2 := \{z \in \Omega_1 : \varphi_j (z) \in \Omega_1, 
                     1 \le j \le n-1 \} . 
\] 
                     Then $\Omega_0 \setminus \Omega_2$ (and therefore 
$\Dd \setminus \Omega_2$) only has isolated points. 

Indeed, $z \in \Omega_0 \setminus \Omega_2$, if and 
only if there is some  $j \in \{0,\dots,n-1\}$ such that $w=\varphi_j(z)$, 
$w \in \Omega_0 \setminus \Omega_1$. Take a point $z_0 \in \Omega_0$, 
and a compact neighborhood of $U$ of $z_0$; then $B^{-1}(B(U))$ is 
a compact set (since a finite Blasche product is a proper map from the 
disk to itself), therefore can only contain a finite number of points 
from $\Omega_0 \setminus \Omega_1$ by the above claim, so there can 
only be a finite number of points in $\Omega_0 \setminus \Omega_2$ in 
$U$, q.e.d. 

                     For any point $z \in \Omega_2$, the fiber 
$B^{-1}(B(z))$ is made up of 
                     exactly $n$ points, and each point in this fiber is 
part of a group of 
                     $m$ points for which the values of $g$ coincide. So 
the theorem holds with 
                     $\Omega=\Omega_2$. 
                     \end{proof}

\section{$z$-invariant subspaces in the algebra's closure}
\label{sec crit}
                   
In general one may ask, first, how large can be the 
$z$-invariant subspaces contained in the closure of $\cala=\Cc [B,g]$? and, second, can 
$\Cc [B,g]$  be included in some closed space that is easy to describe in 
terms of the well-understood $\Cc [B]$ and the function $\Gamma(g)$?  In the introduction 
both questions were answered in the case of Blaschke products of degree two. In 
this section we consider the case of Blaschke products of higher degree.

As to the second question, all we can say for now is 
that $\Gamma(g)$ can be expressed as a polynomial of degree at most 
$n-1$ in $g$ with coefficients in $\overline{\Cc[B]}$ (this follows from the 
proof of Theorem \ref{main claim}, see below). Therefore all elements of $\Cc [B,g]$ are algebraic 
of degree at most $n-1$ over the ring 
$\overline{\Cc[B]}[\Gamma(g)]\subset \overline{\Cc[B]}+ \overline{\Gamma(g)H^p}$. 

The first question is not easy to answer in full 
generality.  The following examples, reflecting respectively cases where the inner factor 
of $\Gamma(g)$ has either zeros inside the disk and singular factors at the 
boundary, suggest that $\Gamma (g) H^p$ is optimal in certain 
cases. 

\vskip5mm 

{\bf Example 1.}  Suppose $B(z)=z^n$ and $g(z)=z^m$. 

Then if the greatest common divisor of $n$ and $m$ is 
$\delta \ge 2$, both $B$ and $g$ are expressed in terms of 
the function $u(z)= z^\delta$, $\Gamma (g)(z) \equiv 
0$ and $\overline{\cala}$ contains no nontrivial invariant subspace. 

So assume that $n$ and $m$ are relatively prime.   
Then $\Gamma (g)(z) = c z^{(n-1)(m-1)}$, 
where $c=c(n,m)$ is a nonzero complex constant. 

Here $\Cc[B,g]$ is the span of the monomials $\{z^k, k \in A\}$, where 
$A:=\{an+bm , a,b \in \bf  Z_+\}$, and its closure is made up of the functions $f\in 
H^p$ such that $f^{(k)}(0)=0$ for $k \notin A$. 

It is easy to see that the largest invariant subspace 
that can be contained in  the closure of this algebra will be of the form 
$z^N H^2$, where $N$ is the smallest integer such 
that any $\{N,N+1,\dots\}\subset A$. 
It is a classical result in arithmetic that $N=(n-1)(m-1)$ (cf \cite[Theorem 
1.17, p. 39]{MBN}, see also \cite {JLB}). 

Note also that in this case the full algebra 
$\Cc[B,g]$ is obtained by adjoining the $(n-1)$-th root of $\Gamma(g)$ to 
$\Cc[B]$ (so in this case the comments made above about the second question 
give the complete answer). The codimension of the closure of $\Cc[B,g]$ is 
$(1/2)(n-1)(m-1)$  (cf \cite[exercise 7, p. 41]{MBN}, or \cite {JLB}).

{\bf Example 2}. Let $B(z)=z^n$ and $g(z)=z \chi 
(z^n)$ (as above $\chi $ is the singular inner function determined by the unit 
point mass at $1$). Then $\Gamma (g)(z) =  \chi (z^n)^{n-1}=S_g(z)$. 
Note that  for any $\al <1/(n-1)$, 
$\Gamma (g)^\al H^p \not \subset \overline{\cala}$. 

Indeed, let $f(z):= z \Gamma (g)(z)^\al = z \chi 
(z^n)^{(n-1)\al}$,  then $\Gamma(f)(z) =  \chi (z^n)^{(n-1)^2\al}$, which 
is not divisible   by $S_g$ since $(n-1)^2\al < n-1$. Thus, by 
Proposition \ref{Aincl}, $f \notin \overline{\cala}$. 

                     \vskip5mm 

 There is a partial necessary condition about the 
vanishing of any inner function $J$ such that $JH^p \subset \overline{\cala}$. 

 \begin{proposition} 
Let $\overline{\cala}=\overline{\Cc[B,g]}$, $B$ a 
Blaschke product of  degree $n$. 
Suppose that $JH^p \subset \overline{\cala}$, where 
$J$ is inner.  Then, for any $a \in \Dd$ such that $\Gamma(g)(a)=0$, 
we have $J(a)=0$.
   \end{proposition} 
\begin{proof} 
 Let $a$ be a zero of $B_g$ (recall that $B_g$ is the 
Blaschke part of the canonical factorization of $\Gamma (g)$).

{\bf Case 1.}
Suppose that $a\neq \varphi_j(a)$, for $1\le j\le 
n-1$. 

Then set $f(z):= \prod_1^{n-1} (z-\varphi_j(a))$. By our 
hypothesis,  $fJ \in \overline{\cala}$, so by Proposition 
\ref{Aincl}, \  $\Gamma(fJ)$ is divisible by $B_g$ . 

We have 
\[ 
\Gamma(fJ) (z) = \prod_{j=1}^{n-1} \frac{f(z)J(z) - f 
\circ \varphi_j (z) J \circ \varphi_j (z)} 
{z - \varphi_j (z)}, 
\] 
and $f \circ \varphi_j (a)=0$ for $1 \le j \le n-1$. 
So 
\[ 
\Gamma(fJ) (a) = \frac{(f(a)J(a))^{n 
-1}}{\prod_{j=1}^{n-1} (a - \varphi_j (a))} =0, 
\] 
which implies $J(a)=0$. 

{\bf Case 2.} Suppose that $a =\varphi_j(a)$ for some $j$.

Now $B_g(a)=0$ is due to  the fact that $B'(a)=g'(a)=0$. Then any $f \in 
\overline{\cala}$ must satisfy $f'(a)=0$, but if we had $J(a)\neq 0$, then for 
$f(z):=z J(z)$ we would have $f'(a)\neq 0$, 
so $J H^p \not \subset \overline{\cala}$. 
\end{proof}

\section{Applications to the disk algebra}
\label{sec cont case}

Let us consider the special case where the function $g$ is a member
of the disk algebra, that is to say $g \in A(\overline \Dd)
:= \mathcal O (\Dd) \cap \mathcal C (\overline \Dd)$, the algebra of 
holomorphic functions continuous up to the boundary of the disk, endowed 
with the uniform norm on $\overline \Dd$. One is now interested 
in $\overline \cala$, where the closure is taken in $A(\overline \Dd)$.

The closure of $\Cc [f,g]$ in $A(\overline \Dd)$ has been studied by 
several authors, see e.g. \cite{JoW}, \cite{RGB}. We can apply
the methods of the previous sections in the case
when $g$ is continuous up to the unit circle. Notice in particular that 
then $\Gamma(g) \in  A(\overline \Dd)$, and  
whenever $\Gamma(g)$ is zero-free, 
$\Gamma(g) A(\overline \Dd) =  A(\overline \Dd)$.

\begin{theorem}
\label{cont case}
$\overline{\cala}$ contains
$\Gamma (g) A(\overline \Dd)$.
\end{theorem}

\begin{corollary}
\label{NaSCcont}
\begin{enumerate}
\item
$\cala$ is dense in $A(\overline \Dd)$ 
if and only if $\Gamma(g)$ has no zeros in the closed unit disk
(in particular $\Gamma(g)$ is then outer).
\item
$\overline {\cala}$ has finite codimension in $A(\overline \Dd)$ if and 
only if 
$\Gamma(g)$ has no 
singular inner factor, and only a finite number of zeros in the 
closed unit disk.
\end{enumerate}
\end{corollary}

 The proofs of both results rely on the proof of
Theorem \ref {main claim} so we will postpone them  until the end of  
section \ref {sec proof}.

\vskip.5cm

{\bf Comparison with previous results}

The earliest results which we know in this direction have been obtained
by John Wermer in \cite{JoW}, where the uniform closure of an algebra
of the form $\Cc [f,g]$ was considered, with $f, g$ analytic in a 
neighborhood of the closed disk (actually, Wermer's results hold in the 
more general context of a domain bounded by a simple curve on a 
Riemann surface) and $f'(\zeta) \neq 0$ for any $\zeta \in \partial \Dd$.

Then Theorem 1.1 in \cite{JoW} states that $\overline{\Cc [f,g]}
=A(\overline \Dd)$ if and only if the algebra separates points 
of the closed disk and the
derivatives of the generators never vanish simultaneously. 

In our case (where $f=B$, a finite Blaschke product), by Lemma
\ref{zeroGamma}, this is equivalent to $\Gamma (g)$ having no zeros
on the closed disk (a condition which automatically rules out singular
inner factors). So a more special hypothesis on the function $f$ 
helps us relax considerably the hypotheses on $g$.

It should be noted that consideration of the function
$g_r(z):= g(rz)$, for $r<1$, which is analytic in a neighborhood
of the closed disk, does not help in this instance, because the
convergence of $g_r$ to $g$ is not enough to ensure the corresponding
``convergence" of the closed algebras they generate. 

For instance, consider the following example:
\[
B(z)=z^2, \quad
g(z)  = z (1 - z^2) \chi(z^2) ,
\]
where $\chi$ stands for the same singular inner function as in the 
examples of the introduction. Then $\Gamma (g)(z)= (1 - z^2) \chi(z^2)$,
so it has a singular inner factor and $\overline{\Cc[B,g]}$ is of 
infinite codimension. On the other hand, 
$ \Gamma (g_r)(z)=r(1 - r^2 z^2) \chi(r^2 z^2)$ has no zeros on
the closed disk for any $r<1$, so then 
$\overline{\Cc[B,g_r]}= A(\overline \Dd)$.

In \cite{RGB}, more general algebras $\cala$ are considered. If 
$\cala$ separates points of the closed disk,
contains a dense subalgebra of functions which are 
${\mathcal C}^1$ up to the boundary, and for any $z \in \overline{\Dd}$
there is $f_z \in \cala$ such that $f'_z (z) \neq 0$, then 
$\cala$ is dense in $A(\overline {\Dd})$. Actually, this immediate
consequence of Blumenthal's work is quoted thus in \cite{JWe}, a paper 
which summarizes and slightly generalizes
an unpublished dissertation \cite{JJo}. As far as we know, those were
the first works to point out the role as obstructions to density of 
singular inner factors of functions analogous to our $\Gamma(g)$.

For algebras of the special form $\overline{\Cc[B,g]}$, Blumenthal's
result extends Wermer's theorem to the case where 
$g \in {\mathcal C}^1(\overline \Dd)$. In this case, our methods can cover
further examples of functions $g$ which yield a dense
subalgebra of $A(\overline {\Dd})$.  For instance, consider
\[
B(z)=z^2, \quad
g(z)=(1-z) \chi(z) \in A(\overline {\Dd}) \setminus 
{\mathcal C}^1(\overline \Dd).
\]
Then 
\[
\Gamma(g)(z) = \frac1{2z} \left( (1-z) \chi(z) - (1+z) \chi(-z) \right).
\]
One checks easily that $\Gamma(g)(0) \neq 0$. 

In fact, $\Gamma(g)$ never vanishes. For $z \neq 0$,
setting $\sigma:= (1+z)/(1-z)$, we see that $\Gamma(g)(z) = 0$
if and only if 
\[
\exp (\frac1\sigma - \sigma) = \sigma,
\]
and since Re$(\frac1\sigma - \sigma) > 0$ if and only if
$|\sigma|<1$, any solutions would have to be of the form 
$\sigma = e^{i\theta}$, which leads to $e^{-2i\sin\theta} = e^{i\theta}$,
an equation whose only solution in the relevant interval 
$(-\frac\pi2,\frac\pi2)$ is $\theta=0$, which corresponds to
$z=0$, excluded.

Furthermore, $\Gamma(g)$ is analytic in a neighborhood of 
$\overline \Dd \setminus \{+1, -1\}$, so any singular inner factor 
would have its singular mass carried by $\{+1, -1\}$. However, one
verifies immediately that $\Gamma(g)$ has nonzero limits at each
of those points, so that it cannot admit any singular factor.
Therefore, by Corollary \ref{NaSCcont}(i),
 $\overline{\Cc[B,g]} = A(\overline {\Dd}) $.
 
\vskip.5cm
The papers of Wermer and Blumenthal cited above also 
include results of the same type as Theorem \ref{cont case},
whose conclusion is that the closure of a given algebra contains
a certain ideal of finite codimension. 

Theorem 1.2 in \cite{JoW} states 
that, under the same regularity hypotheses on $f,g$,
if $\Cc [f,g]$ separates points on $\partial \Dd$, then
$\overline{\Cc [f,g]}$ contains a subspace of the form 
$P A(\overline {\Dd})$, where $P$ is a polynomial.
It should be noted that the hypothesis about point-separation
is easily weakened here to the hypothesis that the set of points
which are not separated be discrete; then a small change in the radius of 
the
boundary circle (legitimate since our functions are analytic in
a neighborhood of the closed disk) will produce a situation where
the points of the circle are separated.

In the case of $\Cc [B,g]$, if $g$ is analytic in a neighborhood
of the closed disk, then so is $\Gamma(g)$, so it can't have 
singular inner factors, nor an infinite number of zeros on the
closed disk unless it is identically zero. That last occurrence 
rules out any non-trivial closed ideal being contained in $\Cc [B,g]$,
by the proof of Corollary \ref{vanish ident}.

 Blumenthal's Theorem 1 in \cite{RGB} extends Wermer's theorem
to cases where $g \in {\mathcal C}^1 (\overline \Dd)$ only; in view of 
Theorem \ref{cont case} this assumption is
not needed when considering $\Cc [B,g]$. (Furthermore, Blumenthal's 
theorem 
includes a hypothesis of point-separation which is 
stronger than expected for the conclusion obtained).

An outstanding problem is the determination of 
$\overline{\Cc [f,g]}$ when $f$ and $g \in A(\overline {\Dd})$,
even if $f$ is assumed to be analytic in a neighborhood of $\overline 
{\Dd}$
(but no extra hypotheses are imposed on $g$).

\section{Proof of the Main Theorem}
\label{sec proof}
\begin{proof}[Proof of Theorem \ref{main claim}] 
Suppose $f \in H^p$, let $f_r (z) = f(rz)$ for $r<1$. 
If a sequence $f_\nu$ can be found in $\cala$ 
which approximates $f_r$ pointwise boundedly, 
then $f_r$ belongs
to the weak closure of $\cala$ in $H^p$. Since the weak closure of $\cala$
coincides with the strong closure, $f_r$ is in the $H^p$-closure of $\cala$. 
 The proof of Theorem \ref{main claim} therefore 
reduces to the following Proposition. 
                     \end{proof} 

                     \begin{proposition} 
\label{main step}
                     For any $h \in \Gamma (g) H^\infty$, there is a 
sequence $\{h_\nu\} \subset \cala$ converging to 
                     $h$ pointwise boundedly. 
                     \end{proposition} 

 \begin{proof} 
First we need to introduce notations for two kinds of 
singularities (in a certain sense) that we will encounter. 
Recall from (\ref{defso}) that $ \cals_0$ is
 the (finite) set of $z$ such that $\{ \varphi_j (z), 0\le j \le n-1\}$ is made up of 
fewer than $n$ points. 

Larger still is the set $\cals_1$ of $z$ such that 
$\{ g(\varphi_j (z)), 0\le j \le n-1\}$ is made up of 
fewer than $n$ points. The previous considerations 
show that $\cals_1 \setminus \cals_0$ is contained in the zero 
set of $\Gamma(g)$, so that if it is not discrete, 
$\Gamma(g)\equiv 0$,  and our theorem is vacuously true.

Let $h$ be any holomorphic function in the disc (later we 
will  choose $h= f \Gamma(g) $).

Restrict attention for a moment to points $z \in \Dd 
\setminus  \cals_1$. Then there is a unique polynomial 
$\call $ of degree not exceeding $n-1$ such that 
\[ 
\call (g \circ \varphi_j (z)) = h(\varphi_j (z))  , 
0\le j \le n-1 . 
\] 
                     The polynomial $\call$ can be written 
                     $\call(X) = \sum_{k=0}^{n-1} \al_j (B(z)) X^k$, since 
its coefficients 
                     depend only on the set $B^{-1}(B(z))$, and, being 
symmetric in the 
                     elements of that set and given by fractions 
(Lagrange's interpolation 
                     formula), are holomorphic in $\Dd \setminus \cals_1$. 
So we have, 
                     at any such point $z$, 
                     \[ 
                     h(z) = \sum_{j=0}^{n-1} \al_j (B(z)) g(z)^j . 
                     \] 

 Now it suffices to show that when $h=f \Gamma (g)$, 
the $\al_j$ turn out to be bounded analytic 
functions in the unit disk.   We write $\call = \call_z$ to emphasize the 
dependency on $z$. 

 In a punctured neighborhood of any point of 
$\cals_1$, $\call$ is given by the Lagrange formula, 
\begin{equation} 
\label{Lagrange} 
\call_z (X) = 
\sum_{j=0}^{n-1} f(\varphi_j(z)) \Gamma (g) (\varphi_j (z)) 
\prod_{l:l\neq j} \frac{X- g(\varphi_l (z))} { g(\varphi_j (z))- g(\varphi_l (z))} . 
\end{equation} 
On the other hand, 
\[ 
\Gamma (g) (\varphi_j (z)) = 
\prod_{i=1}^{n-1} \frac{ g(\varphi_j (z))- 
g(\varphi_i \circ \varphi_j (z))} {\varphi_j (z)- \varphi_i \circ \varphi_j (z)} . 
\] 
When $i$ goes through the set $\{1,\dots,n-1\}$, 
$\varphi_i \circ \varphi_j (z)$ goes through the set 
$\{ \varphi_l (z), l \neq j\}$, so the above numerators cancel out 
with the denominators  (\ref{Lagrange}). 
 So finally 
\begin{equation} 
\label{LagrSimp} 
\call_z (X) = 
\sum_{j=0}^{n-1}  f(\varphi_j (z)) 
\prod_{l:l\neq j} \frac{X- g(\varphi_l (z))} 
{\varphi_j (z)- \varphi_l (z)} . 
\end{equation} 
This formula makes  sense whenever $z \notin \cals_0$, and shows that the 
coefficients of $\call$ are holomorphically extendable to $\Dd 
\setminus \cals_0$. 

     Consider now a neighborhood $V$ of the unit circle 
chosen so that its closure is disjoint from 
$\cals_0$.  For $z \in V$,  $|\varphi_j (z)- \varphi_k (z)|$ is bounded below. 
So to prove that the $\alpha_j$ are bounded analytic 
functions,  it will be enough to prove that the singularities 
for $z \in \cals_0$ are removable. 

This will be done by considering several cases in 
turn for the possible singularity $z_0 \in \cals_0$. 
First, we consider   a type of function $g$ which allows us to describe 
explicitly  what the interpolating polynomial $\call_z$ turns out 
to be when  $z$ is a critical point of $B$. 

                     \vskip5mm 

                     {\bf Case 1.}  $g'(z) \neq 0$ for any $z \in 
B^{-1}(B(z_0))$, and 
                     $g$ is one-to-one on the set $B^{-1}(B(z_0))$. 

                     Note that this hypothesis is equivalent to the fact 
that 
                     $g$ is one-to-one in some neighborhood $U$ of 
$B^{-1}(B(z_0))$. 

                     Then $g^{-1}$ is holomorphic on the open set $g(U)$, 
and $\call$ 
                     is exactly the interpolation polynomial of degree 
less or 
                     equal than $n-1$ assuming the same values as $h \circ 
g^{-1}$ 
                     at the points of $g \circ B^{-1}(B(z))$, for $z \in U 
\setminus \cals_0$. 
                     Choose a smaller neighborhood $U_0$ of 
$B^{-1}(B(z_0))$, 
                     $U_0$ relatively compact in $U$, such that the 
boundary of $g(U_0)$ is 
                     a finite disjoint union of circles. For $z$ belonging 
to a small enough 
                     neighborhood of $z_0$ such that $B^{-1}(B(z)) \subset 
U_0$, 
                     and for $x \in B(U_0)$, 
\[ 
                     \call_z (x) = 
                     \frac1{2\pi i} \int_{\partial U_0} 
                     \frac{h \circ g^{-1} (\zeta)}{\omega_z (\zeta)} 
                     \frac{\omega_z (\zeta) - \omega_z (x)}{\zeta - x} \,d 
\zeta , 
\] 
where $\omega_z (X) := \prod_{0\le l \le n-1} \left( 
X- g\circ\varphi_l (z) \right)$ \cite[Vol. II, Section 11, p. 67 ff.]{AIM}. This formula 
depends continuously on $z$ in the whole neighborhood 
$U_0$, and converges, when $z$ converges to $z_0$, 
to the solution of an Hermite interpolation 
problem given by 
\[ 
\call^{(j)} (g \circ \varphi_i (z)) = h^{(j)} 
(\varphi_i (z)), 
{\rm\ for\ }i \in A_q, 0\le j \le \# A_q -1, 1 \le q 
\le p, 
 \] 
where 
$A_1, \ldots, A_p$ form a partition of 
$\{0,\ldots,n-1\}$ such that $\varphi_j (z) = 
\varphi_{j'} (z)$ 
if and only if $j$ and $j'$ belong to the same set 
$A_q$. 

                     \vskip5mm 

{\bf Case 2.}  For any point $z_0$ such that there 
exists a $j$ with $B'(\varphi_j(z_0))=0$, we have 
 $B''(\varphi_j(z_0)) \neq 0$ for any such $j$. 

This hypothesis says that the equation $B(z)=B(z_0)$ 
has solutions  of multiplicity never exceeding $2$. It does not 
depend on $g$.  Notice that this case is generic: Blaschke products 
 of degree greater than one always admit critical 
points inside the disk, and those for which 
multiplicities of the above solutions can be $3$ or more form a small set,
 in a sense that will be made precise at the beginning of Case 3. 

Then the various terms in the sum which represents 
$\call_z(X)$ either converge to a finite limit (when $j$ is such that 
$B'(\varphi_j(z_0))\neq 0$) or can be regrouped in pairs $(j,k)$ where $j \neq k$, 
$\varphi_j(z_0) = \varphi_k(z_0) \neq \varphi_l(z_0)$ for $l \notin \{j,k\}$. To 
simplify notations, let us assume that $j=0$, $k=1$. 
For $z$ in a small enough neighborhood of $z_0$, the sum of the two terms can be written 
\begin{eqnarray} 
\label{twotermsum}  
&&f(z) \frac{X-g \circ \varphi_1(z)}{z-\varphi_1(z)} 
\prod_{2 \le l \le n-1 } \frac{X- g \circ \varphi_l (z)} {z - \varphi_l (z)} 
                     \\ 
\nonumber 
&& \quad 
+ 
f\circ\varphi_1(z) 
\frac{X-g (z)}{\varphi_1(z)-z} 
\prod_{2 \le l \le n-1 } \frac{X- g \circ \varphi_l (z)} {\varphi_1(z) -  \varphi_l (z)} 
\\ 
\nonumber
&=& 
\frac{\prod_{2 \le l \le n-1 } (X- g \circ \varphi_l 
(z))}{z- \varphi_1 (z)} 
\left( 
\frac{f(z)(X- g \circ \varphi_1 (z))}{\prod_{2 \le l \le n-1 } (z- \varphi_l (z))} 
- 
\frac{f\circ \varphi_1 (z) (X- g (z))} 
{\prod_{2 \le l \le n-1 } ( \varphi_1 (z)-  \varphi_l (z))} 
\right) 
. 
\end{eqnarray} 
Define a meromorphic function whose coefficients depend on $z$ by 
\[ 
R(Y):= \frac{f(Y)}{\prod_{2 \le l \le n-1 } (Y- \varphi_l (z))}. 
\] 
 The poles of $R$ depend on the parameter $z$ and stay away from $z_0$ for $z$ in a 
neighborhood of $z_0$, and we have 
\[ 
R(Y_1)-R(Y_2):= (Y_1-Y_2) R_1 (Y_1,Y_2) , 
\] 
where $R_1$ is another meromorphic function (of two 
variables) with its pole set (which is a hypersurface of 
$\mathbb C^2$) avoiding a neighborhood of $(z_0,z_0)$.
 The expression from
(\ref{twotermsum}) is then equal to 
\begin{eqnarray*} 
&&\frac{\prod_{2 \le l \le n-1 } (X- g \circ 
\varphi_l (z))}{z- \varphi_1 (z)} 
\left[ 
R(z) (X- g \circ \varphi_1 (z)) 
- 
R(\varphi_1 (z)) (X- g (z)) 
\right] 
\\ 
&=& 
\left( 
 \prod_{2 \le l \le n-1 } (X- g \circ \varphi_l (z)) 
\right) 
 \left( 
 R_1(z, \varphi_1 (z)) (X- g \circ \varphi_1 (z)) 
+ \frac{g(z)- g \circ \varphi_1 (z)}{z- \varphi_1 (z)} 
R( \varphi_1 (z)) 
\right) . 
\end{eqnarray*} 
Letting $z$ tend to $z_0$, $\varphi_1(z)$ tends to 
$z_0$ as well, and $\frac{g(z)- g \circ \varphi_1 (z)}{z - \varphi_1 
(z)}$ tends to $g'(z_0)$, so the coefficients of $\call_z$ 
admit  finite limits, and the possible singularity of its 
coefficients at  $z_0$ is actually removable. 

                     \vskip5mm 

                     {\bf Case 3.} 
                     General case. 

The Blaschke products of degree $n$ which do not 
satisfy the hypotheses of Case 2 form an algebraic subset, therefore of empty 
interior, since their derivative, which is a rational fraction of fixed 
degree, must admit at  least one multiple root. 

Given one such exceptional Blaschke product $B$, let 
$B^\nu$ be a sequence of Blaschke products satisfying the hypotheses of 
Case 2 and converging to 
$B$ uniformly in a neighborhood 
the closed unit disk.  By continuous dependency of 
the roots  of a polynomial equation upon the coefficients (Hurwitz Theorem), we 
see that the corresponding functions $\varphi^\nu_j$ converge to $\varphi_j$ 
uniformly in a some neighborhood of the unit circle (where they are all well-defined), 
and that any symmetric polynomial in the $\varphi^\nu_j$'s converges 
uniformly in the closed disk to the same polynomial in the $\varphi_j$'s.  In 
particular, if we set 
$\Gamma^\nu(g)=\Gamma_{B^\nu }(g)$,  then 
$\Gamma^\nu(g)$ converges to $\Gamma(g)$ 
pointwise boundedly on the disk. 

Given a function $f \in H^\infty(\Dd)$, write  
\[ 
\call_z^\nu (X) = 
\sum_{j=0}^{n-1} f(\varphi_j^\nu (z)) \Gamma^\nu (g) (\varphi_j^\nu (z)) 
\prod_{l:l\neq j} \frac{X- g(\varphi_l^\nu (z))} 
{ g(\varphi_j^\nu (z))- g(\varphi_l^\nu (z))} =: 
\sum_{k=0}^{n-1} A_k^\nu (z) X^k 
\] 
From (\ref{LagrSimp}) and the proof in Case 2, 
 we see that the $A_k^\nu$ are 
bounded holomorphic functions and $A_k^\nu (z) = 
\al_k^\nu(B^\nu(z))$. 

                     Fix some degree $k$ between $0$ and $n-1$. On any 
 circle centered at the origin of radius $r$ close 
enough to $1$, $A_k^\nu$ converges uniformly, 
because $f$ is bounded on the circle, all the 
$\varphi^\nu_j$ are well defined and converge 
uniformly, and all the denominators in the expression 
are bounded away from $0$. If $f$ is furthermore bounded, the 
limit is a bounded holomorphic function which we 
denote by $A_k$. 

We need to see that for any $z_1$, $z_2$ such that 
$B(z_1)=B(z_2)$, 
                     then $A_k(z_1)=A_k(z_2)$. In that case, there is a 
$j$ such that 
                     $\varphi_j (z_1)=z_2$, and by the convergence of the 
$\varphi_j^\nu$, 
                     we can select a subsequence (which we will denote the 
same way as the original 
                     sequence) and an index $j_0$ such that $z_2 - 
\varphi_{j_0}^\nu(z_1) \to 0$ 
                     as $\nu \to \infty$. Then 
                     \[ 
                     A_k(z_1)-A_k(z_2) = \lim_{\nu\to\infty} 
(A_k^\nu(z_1)-A_k^\nu(z_2)) 
                     = \lim_{\nu\to\infty} 
(A_k^\nu(\varphi_{j_0}^\nu(z_1))-A_k^\nu(z_2)) =0, 
                     \] 
                     where the next to last equality is by the symmetry of 
the $A_k^\nu$, and the 
                     last by uniform convergence on compacta.  So we can 
write 
                     $A_k (z) = \al_k (B(z))$, where $\al_k \in 
H^\infty(\Dd)$. 
                     Now 
                     \[ 
                     f(z) \Gamma(g)(z) =  \lim_{\nu\to\infty} f(z) 
\Gamma^\nu(g)(z) 
                     =  \lim_{\nu\to\infty} \sum_{k=0}^{n-1} A_k^\nu (z) 
g(z)^k 
                     = \sum_{k=0}^{n-1} \al_k (B(z)) g(z)^k. 
                     \] 
                     \end{proof} 

We now are able to give the proofs of the corresponding results for the disk 
algebra case.

\begin{proof}[Proof of Theorem \ref{cont case}]
The proof proceeds exactly as the proof of Proposition \ref{main step},
except that now we will need to assume $h= f \Gamma(g)$, with 
$f \in A(\overline \Dd)$.
This time we need to prove that the functions $\alpha_j$ are continuous
up to the boundary (in addition to being holomorphic). For this it will be
enough to show that the coeffients of $\call_z$ are themselves continuous
up to the boundary. This follows from the formula (\ref{LagrSimp}), 
since $f$ and $g$ are continuous, and the $\varphi_j$ are all distinct 
near 
the unit circle.

Furthermore, the arguments of Case 3 now show that 
$\Gamma^\nu(g)$ converges to $\Gamma(g)$
uniformly on the closed disk (using the uniform continuity of
$g$ on the closed disk).  So the coefficients $A_k^\nu$ converge
uniformly to functions in the disc algebra, and one proves in the 
usual way that the limits are functions of $B$ alone.
\end{proof}

\begin{proof} [Proof of Corollary \ref{NaSCcont}]

(i) Any zero of $\Gamma(g)$ in the open disk provides a continuous, 
non-trivial
linear form which vanishes on $\cala$, by the proof of Corollary 
\ref{NCGamma}. If $\Gamma(g)$ has a zero $\zeta$ on the boundary, then 
since
$B'$ never vanishes on the unit circle, there must exist 
$\zeta' \in \partial \Dd$, $\zeta' \neq \zeta$, such that 
$B(\zeta)=B(\zeta')$, $g(\zeta)=g(\zeta')$. Then 
$f \mapsto f(\zeta)-f(\zeta')$ provides a non-trivial
linear form, continuous on $A(\overline \Dd)$, which vanishes on $\cala$.

The converse is immediate by the remark made before the statement of
Theorem \ref{cont case}.

(ii) As above, an infinite number of zeros of $\Gamma(g)$ would provide 
arbitrarily large (finite) sets of linearly independent 
continuous linear forms vanishing on $\cala$, so the codimension could not
be finite.

Suppose on the other hand that the factorization of $\Gamma(g)$ admit
a singular inner factor $S$. It is easy to see that the map 
$f \mapsto \Gamma(f)$ is continuous on $A(\overline \Dd)$, so that
whenever $f \in \overline{\cala}$, $\Gamma(f)$ is divisible by $\Gamma(g)$,
therefore by $S$. Arguing as in the proof of Corollary 
\ref{NCGamma} again, if $\overline{\cala}$ was of finite codimension,
it would contain a polynomial, which would lead to a contradiction.

Conversely, if $\#\Gamma(g)^{-1}\{0\}$ is finite, it will be enough to
show that the codimension of the closure of $\Gamma (g) A(\overline \Dd)$
is finite. Let ${\mathcal I} := \overline{\Gamma (g) A(\overline \Dd)}$. 
Then
$\mathcal I$ is a closed ideal of $A(\overline {\Dd})$, and the main result
of \cite {WRu}
shows that ${\mathcal I}= I(E,M)$, the set of all functions in 
$A(\overline {\Dd})$
which vanish on the closed set $E \subset \partial \Dd$ and which are
divisible by the inner function $M$. The hypothesis on $\Gamma(g)$
imply that $M$ is a finite Blaschke product
(because singular inner factors are excluded), and that 
$E= \Gamma(g)^{-1}\{0\}\cap \partial \Dd$ is a finite
subset of the circle. Therefore $\mathcal I$ reduces to those functions
in $A(\overline \Dd)$ which vanish on $E':=\Gamma(g)^{-1}\{0\}\cap \Dd$
to an order greater or equal than the order of vanishing of 
$\Gamma(g)$ at those points, and which vanish on $E$ (where the notion of 
``order of vanishing" makes no sense). 

Notice that this allows us to compute precisely the codimension of
$\mathcal I$, as the sum of the number of points in $E$, counted
in the usual way, and the number of
points in $E'$, counted with multiplicities.
\end{proof}

\vspace{1cm}

$$
\begin{array}{ll}
\mbox{Michael I. Stessin} & \mbox{ \ \ \ Pascal J. Thomas} \\
\mbox{Department of Mathematics and Statistics} & \mbox{ \ \ \ Laboratoire 
Emile Picard, UMR CNRS 5580} \\
\mbox{University at Albany} & \mbox{ \ \ \ 118 route de Narbonne} \\ 
\mbox{Albany, NY 12222} & \mbox{ \ \ \ 31062 TOULOUSE CEDEX} \\ 
\mbox{USA} & \mbox{ \ \ \  France} \\ 
\mbox{stessin@math.albany.edu} & \mbox{ \ \ \ pthomas@cict.fr} \\
\end{array}
$$

                    \end{document}